\documentclass[11pt,twoside]{article}
\title{Rings over which all (finitely generated strongly) Gorenstein projective
 modules are projective}
\date{}
\usepackage{amsfonts}
\usepackage{amsmath}
\usepackage{amssymb}
\usepackage{latexsym}
\usepackage{graphicx}

\newtheorem{thm}{\bf Theorem}[section]
\newtheorem{cor}[thm]{\bf Corollary}
\newtheorem{lem}[thm]{\bf Lemma}
\newtheorem{prop}[thm]{\bf Proposition}

\newtheorem{defns}[thm]{\bf Definitions}

\newtheorem{exmp}[thm]{\bf Example}

\catcode`\ç=13
\defç{\c{c}}
\catcode`\é=13
\defé{\'e}
\catcode`\à=13
\defà{\`a}
\catcode`\è=13
\defè{\`e}
\catcode`\â=13
\defâ{\^a}
\catcode`\ù=13
\defù{\`u}
\catcode`\ê=13
\defê{\^e}
\catcode`\î=13
\defî{\^\i}
\catcode`\ô=13
\defô{\^o}

\def\proof{{\parindent0pt {\bf Proof.\ }}}

\def\wdim{{\rm wdim}}
\def\gldim{{\rm gldim}}

\def\Ggldim{{\rm G\!-\!gldim}}

\def\pd{{\rm pd}}

\def\id{{\rm id}}

\def\Gpd{{\rm Gpd}}

\def\Im{{\rm Im}}

\def\Ext{{\rm Ext}}

\def\Hom{{\rm Hom}}

\newcommand{\cqfd}
{\hspace{1cm}
\rule{2mm}{2mm}%
\medbreak%
\par%
}
\def\1{{\noindent\rm (1)}}
\def\2{{\noindent\rm (2)}}
\def\3{{\noindent\rm (3)}}
\def\4{{\noindent\rm (4)}}
\def\5{{\noindent\rm (5)}}
\begin{document}

\centerline{}

\centerline{}

 \centerline {\Large{\bf Rings over which all (finitely generated strongly) Gorenstein}}
\centerline {\Large{\bf projective modules are projective}}
\thispagestyle{empty} \centerline{}

\centerline{\bf {Najib Mahdou}}

\centerline{}

\centerline{Department of Mathematics, Faculty of Science and
Technology of Fez,}

\centerline{Box 2202, University S. M. Ben Abdellah Fez, Morocco}

\centerline{mahdou@hotmail.com}

\centerline{}

\centerline{\bf {Khalid Ouarghi}}

\centerline{}

\centerline{Department of Mathematics, Faculty of Sciences, King
Khaled University}

\centerline{PO Box 9004, Abha, Saudi Arabia}

\centerline{ouarghi.khalid@hotmail.fr}

\centerline{}

\bigskip\bigskip

\noindent{\large\bf Abstract.} The main aim of this paper is to
investigate rings over which all (finitely generated strongly)
Gorenstein projective modules are  projective.  We consider this
propriety under change of rings, and give various examples of
rings with and without this propriety.\bigskip

\bigskip\bigskip



\begin{section}{Introduction}
Throughout this paper all rings are commutative with identity
element and all modules are unital. For an $R$-module $M$, we use
$\pd_R(M)$ to denote the usual projective dimension of $M$.
$\gldim(R)$ and $\wdim(R)$ are, respectively, the classical global
and weak global dimensions of $R$. It is convenient to use
``$m$-local" to refer to (not necessarily Noetherian) rings with a
unique maximal ideal $m$.

In 1967-69, Auslander and Bridger \cite{A1,A2} introduced the
  G-dimension for finitely generated modules over Noetherian rings.
Several decades later, this homological dimension was extended, by
Enochs and Jenda \cite{GoInPj,GoIn}, to the Gorenstein projective
dimension of modules that are not necessarily finitely generated
and over non-necessarily Noetherian rings. And, dually, they
defined the Gorenstein injective dimension. Then, to complete the
analogy with the classical homological dimension, Enochs,  Jenda
and  Torrecillas \cite{GoPlat} introduced the Gorenstein flat
dimension.

In the last years, the Gorenstein homological dimensions have
become a vigorously active area of research (see \cite{LW,Rel-hom}
for more details). In 2004, Holm \cite{HH} generalized several
results which already obtained over Noetherian rings. Recently, in
\cite{BM} the authors introduced  particular cases of Gorenstein
projective, injective, and flat modules, which are called
respectively, strongly Gorenstein projective, injective and flat
modules, which are defined, respectively, as follows:

\begin{defns}[\cite{BM}]\label{defSG}
\begin{enumerate}
    \item  A module $M$ is said to be \textit{strongly
Gorenstein projective}, if there exists a complete projective
resolution of the form $$ \mathbf{P}=\
\cdots\stackrel{f}{\longrightarrow}P\stackrel{f}{\longrightarrow}P\stackrel{f}{\longrightarrow}P
\stackrel{f}{\longrightarrow}\cdots $$ such that  $M \cong
\Im(f)$.
   \item The strongly Gorenstein injective modules are defined dually.
    \item A module $M$ is said to be \textit{strongly
Gorenstein flat}, if there exists a complete flat resolution of
the form $$ \mathbf{F}=\
\cdots\stackrel{f}{\longrightarrow}F\stackrel{f}{\longrightarrow}F\stackrel{f}{\longrightarrow}F
\stackrel{f}{\longrightarrow}\cdots $$ such that  $M \cong
\Im(f)$.
\end{enumerate}
\end{defns}

The principal role of the strongly Gorenstein projective  modules
is to give a simple characterization of Gorenstein projective
 modules, as follows:

\begin{thm}[\cite{BM}, Theorem 2.7]\label{thm-car-G-SG}
A module is Gorenstein projective   if and only if it is a direct
summand of a strongly Gorenstein projective module.
\end{thm}

The important of this last result manifests in showing that  the
strongly Gorenstein projective  modules have simpler
characterizations than their Gorenstein correspondent modules. For
instance:

\begin{prop}[\cite{BM}, Proposition  2.9]\label{pro-cara-SG-pro}
A module $M$ is strongly Gorenstein projective if and only if
there exists a short exact sequence of modules: $$0\rightarrow
M\rightarrow P\rightarrow M\rightarrow 0,$$ where $P$ is
projective, and $\Ext(M,Q)=0$ for any   projective module $Q$.
\end{prop}

In order to give an answer to the question ``when is a finitely
generated torsionless module projective?" Luo and Huang  proved,
for a commutative Artinian ring $R$, that a Gorenstein projective
$R$-module $M$ is projective if $\Ext_R^i(M,M)=0$ for any $i\geq
1$ \cite[Theorem 4.7]{LH}. And over Noetherian local rings
Takahashi proved that $G$-regular rings are the rings over which
all $Gdim(M)=\pd(M)$ for any $R$-module $M$ \cite[Proposition
1.8]{Tak}. Over a  local ring $(R,m)$  satisfies that $m^2=0$ and
$R$ is not a Gorenstein ring, Yoshino proved that every $R$-module
of $G$-dimension zero is free \cite[Proposition 2.4]{Yo}. In this
paper, we are concerned with a global question. Namely, we study
the following two classes of rings: rings over which all
Gorenstein projective modules are projective and rings over which
all finitely generated strongly Gorenstein projective modules are
projective.  In Section 2, we show, that first class coincides
with the class of rings over which all strongly Gorenstein
projective modules are projective (see Theorem \ref{pro1}).
Furthermore, in the same result, we show that a ring $R$ belongs
in this class  if and only if $\Ext_R^1(M,M)=0$ for any strongly
Gorenstein projective $R$-module $M$ if and only if
$\Gpd_R(M)=\pd_R(M)$ for any $R$-module $M$.  After, we study  the
second  class over which all finitely generated strongly
Gorenstein projective modules are projective. Then, we study the
transfer of this property in some extensions of rings. In section
3,  we give  some examples of rings with  and without this
property.
\end{section}
\begin{section}{Rings over which all (finitely generated strongly) Gorenstein projective modules are projective}
We start this section with  investigating   rings  satisfy the
property  ``all Gorenstein projective $R$-modules are projective".
In the next Theorem, we see some  conditions equivalent to this
property.

\begin{thm}\label{pro1}
Let  $R $ be a ring. The following conditions are equivalent:
\begin{enumerate}
    \item All Gorenstein projective $R$-modules are projective;
    \item All strongly Gorenstein projective $R$-modules are projective;
    \item For any strongly Gorenstein projective $R$-module $M$,
    $\Ext_R^1(M,M)=0$;
    \item For any $R$-module $M$, $\Gpd_R(M)=\pd_R(M)$.
\end{enumerate}
\end{thm}
\proof  $(2) \Rightarrow (3).$ Is obvious\\
 $(3)
\Rightarrow (2).$ Let $M$ be a  strongly Gorenstein projective
module, from Proposition \ref{pro-cara-SG-pro}, there is an exact
sequence of $R$-modules:
$$0\longrightarrow M\longrightarrow P\longrightarrow M\longrightarrow 0\qquad\qquad (\star)$$
where $P$ is projective. And since $\Ext_R^1(M,M)=0$, the sequence
 $(\star)$
spilt and  $M$
is  a direct summand of $P$, then $M$ is projective.\\
 $(2)
\Rightarrow (1).$ Follows immediately from Theorem
\ref{thm-car-G-SG}.\\
$(1) \Rightarrow (2).$ Obvious (since all strongly
Gorenstein projective modules are Gorenstein projective).\\
$(4)\Rightarrow (1).$ Obvious.\\
$(1)\Rightarrow (4).$ Let $M$ be an $R$-module, it is known that
$\Gpd_R(M)\leq \pd_R(M)$. Then it remains to prove that
$\pd_R(M)\leq\Gpd_R(M)$. If $\Gpd_R(M)=\infty$ it is obvious. Let
$m$ be  a positif integer and $\Gpd_R(M)=m<\infty$. From
\cite[Definition 2.8]{HH}, $M$ has a Gorenstein projective
resolution of length $m$. Then,  $\pd_R(M)\leq m=\Gpd_R(M)$ since
all Gorenstein projective modules are projective. Therefore,
$\Gpd_R(M)=\pd_R(M)$.\cqfd\bigskip

Throughout the remainder of this paper, we study rings satisfy
each of the following conditions equivalent:

\begin{thm}\label{prop 3.1}
Let  $R $ be a ring. The following conditions are equivalent:
\begin{enumerate}
\item All finitely
    generated  strongly Gorenstein projective $R$-modules are
    projective;
\item All finitely generated strongly Gorenstein projective $R$-modules are
    flat;
\item All finitely presented strongly Gorenstein flat $R$-modules are
    projective;
\item All finitely presented strongly Gorenstein flat $R$-modules
    are flat;
\item For any finitely generated strongly Gorenstein projective $R$-module $M$,\\ $\Ext_R(M,M)=0$.
\end{enumerate}
\end{thm}
\proof $(1)\Leftrightarrow (2)$ and $(4)\Rightarrow (1)$. Follows
immediately from \cite[Proposition 3.9]{BM}, and since
every finitely presented flat $R$-module is projective.\\
$(1)\Rightarrow (3).$ Follows from \cite[Proposition 3.9]{BM}.\\
$(3)\Rightarrow (4).$ Obvious.\\
$(1)\Longleftrightarrow (5).$ Similar to the proof of Theorem
\ref{pro1}.\cqfd

Recall,  for an extension of rings $A \subseteq B$, that $A$ is
called a module retract of $B$ if there exists an $A$-module
homomorphism $f : B \longrightarrow A$ such that $f_{/A} =
\id_{/A}$. The homomorphism $f$ is called a module retraction map.
If such map $f$ exists, $B$ contains $A$ as a  direct summand
$A$-module.  In the next  main result, we  study  the property
``all finitely generated strongly Gorenstein projective modules
are projective"
   in   retract rings.

\begin{thm}\label{prop4}
Let  $A$ be a retract subring  of $R$, ($R=A\oplus_A E$), such
that  $E$ is a flat $A$-module. Then, if  the property,  all
finitely generated  strongly Gorenstein projective modules are
projective, holds in $R$, then, it holds in $A$ too.
\end{thm}
\proof   We prove first that $M\otimes_A R$ is a finitely
generated strongly Gorenstein projective $R$-module, for any
finitely generated strongly Gorenstein projective $A$-module $M$.
From \cite[Proposition 2.12]{BM}, there exists an exact sequence
of $A$-modules:
$$0\longrightarrow M\longrightarrow P\longrightarrow M\longrightarrow 0$$
where $P$ is a finitely generated projective $A$-module.
$P\otimes_A R$ is a finitely generated projective $R$-module, and
since $R$ is a flat $A$-module, it follows that the sequence of
$R$-modules:
$$0\longrightarrow M\otimes_A R \longrightarrow P\otimes_A R
\longrightarrow M\otimes_A R \longrightarrow 0$$
 is exact. It remains only to show that $\Ext_R(M\otimes_A R, R)=0$.
Therefore, $\Ext_A(M, R)=0$ ( since $R$ is an $A$-module flat  and
from \cite[Proposition 2.12]{BM}). On the other hand,  from
\cite[Proposition 4.1.3]{CE}, $\Ext_R(M\otimes_A R, R)\cong
\Ext_A(M, R)=0$. Then
 $M\otimes_A R$ is
a finitely generated strongly Gorenstein projective $R$-module and
by hypothesis it is projective. To complete the proof, we will
show that $\Ext_A^k(M,N)=0$ for any integer $k\geq1$ and for any
$A$-module $N$. It is known that  $\Ext_A^k(M,N\otimes_A R)\cong
\Ext_R^k(M\otimes_A R,N\otimes_A R)=0$, (from \cite[Proposition
4.1.3]{CE}).  Namely, $\Ext_A^k(M,N)$ is a direct summand of
$\Ext_A^k(M,N\otimes_A R)$, as $A$-modules (since $A$ is a direct
summand of $R$ as $A$-module). Then,  $\Ext_A^k(M,N)=0$ and $M$ is
a projective $A$-module as desired.\cqfd

Next we study the transfer of the property, all finitely
    generated  strongly Gorenstein projective modules are
    projective, in
polynomial rings.

\begin{cor}\label{thm4}
Let  $R$ be a ring and $X$ an indeterminate over $R$. If $R[X]$
satisfies  the conditions equivalent of Theorem \ref{prop 3.1},
then $R$ satisfies it too.
\end{cor}
\proof Note first  that this Corollary is a particular case of
Theorem \ref{prop4} above, but here we get an other proof. Let $M$
be a finitely generated strongly Gorenstein projective $R$-module.
Since $\pd_R(R[X])$ is finite,  and from \cite[Theorem
2.11]{BMn-SG}, $M[X]$ is a finitely generated strongly Gorenstein
projective $R[X]$-module, so $M[X]$ is a projective $R[X]$-module.
Then, \cite[Lemma 9.27]{Rot} gives that $M$ is a projective
$R$-module.\cqfd

Recall,  let $A$ be a ring and let $E$ an $A$-module. The trivial
ring extension of  $A$ by $E$ is the ring $R := A \propto E$ whose
underlying group is $A \times E$ with multiplication given by
$(a,e)(a',e') = (aa',ae'+a'e)$. We define similarly $J:=I\propto
E'$, where $I$ is an ideal of $A$ and $E'$ is an $A$-submodule of
$E$ such that $IE \subseteq E'$. Then $J$ is an ideal of $R$ and,
if $J$ is a finitely generated ideal, then so is $I$ \cite[Theorem
25.1]{Hu}. Trivial ring extensions have been studied extensively;
the work is summarized in \cite{Triv,Glaz,Hu}. These extensions
have been useful for solving many open problems and conjectures in
both commutative and non-commutative ring theory. See for
instance,\cite{Glaz,Hu,Mh1,Mh2}. As a direct application of
Theorem \ref{prop4} above we have:

\begin{cor}\label{corofprop4}
Let  $A$ be a ring and $E$ a flat $A$-module. If the property, all
finitely generated  strongly Gorenstein projective modules are
projective, holds in $R=A\propto E$,  then
 it holds in $A$ too.
\end{cor}

It is well-known that the structures of ideals are simple more
then the structures of modules, and the study of any property
under ideals gives a large class of examples and solving some
problems. In the following result we study the transfer of the
property ``all finitely generated strongly Gorenstein projective
ideals are projective" between an integral domain $D$ and its
trivial ring extension $D\propto K$ where $K=qf(D)$.

\begin{thm}\label{thmDK}
Let $(D,m)$ be an integral domain $m$-local not field and
$K=qf(D)$. Let $R=D\propto K$, then the following conditions are
equivalent:
\begin{enumerate}
    \item  $R$ satisfies  all finitely generated  strongly Gorenstein projective ideals are projective.
    \item $D$ satisfies all finitely generated  strongly Gorenstein projective ideals are projective.
\end{enumerate}
\end{thm}
To prove this theorem we need the following Lemma.
\begin{lem}\label{lem3.13}
Let $(D,m)$ be an integral domain $m$-local not field and
$K=qf(D)$. Let $R=D\propto K$. If the property,  all finitely
generated strongly Gorenstein projective ideals are projective,
holds in $D$, then $0\propto I$ can not be a strongly Gorenstein
projective ideal of $R$, for any finitely generated ideal $I$ of
$D$.
\end{lem}
\proof Assume, on the contrary, that $ 0\propto I$ is a finitely
generated strongly Gorenstein projective ideal of $R$, for some
finitely generated ideal $I$ of $D$. We prove first   that
$0\propto I$ is  a strongly Gorenstein flat $D$-module. From
\cite[Proposition 2.12]{BM}, there exists an exact sequence of
$R$-modules:
$$0\longrightarrow 0\propto I \longrightarrow P
\longrightarrow 0\propto I\longrightarrow 0\qquad (\star)$$ where
$P$ is finitely generated projective $R$-module, then $P$ is a
free $R$-module (since $R$ is a local ring). Therefore, $(\star)$
is also an exact sequence of $D$-module and since $R$ is a flat
$D$-module,  $P$ is also a flat $D$-module. From \cite[Proposition
3.6]{BM}, it remains to prove that $Tor_D(E, 0\propto I)=0$, for
any injective $D$-module $E$. Thus,  $Tor_D(E, 0\propto I)\cong
Tor_R( \Hom_D(R,E),0\propto I)=0$ (from \cite[Proposition
4.1.1]{CE}). Then, $0\propto I$ is a finitely generated strongly
Gorenstein flat $D$-module and since $D$ is an $m$-local domain
$0\propto I$ is a strongly Gorenstein projective $D$-module (from
\cite[Corollary 3.10]{BM}). On the other hand, $0\propto I\cong I$
as $D$-module. It follows that $I$ is a finitely generated
strongly Gorenstein projective ideal of $D$ and by hypothesis
projective and since $D$ is $m$-local $I$ is  free. Then, $I$ is a
principal ideal of $D$. Let $I=Da$ where $a\in I$. There exists an
exact sequence of $R$-module:
$$0\longrightarrow 0\propto K \longrightarrow R \stackrel{u}{\longrightarrow} 0\propto I\longrightarrow 0 \qquad\qquad (\star\star)$$
where $u((b,e)=(0,ba)$. Therefore, from Schanuel's lemma applied
to the sequences $(\star)$ and $(\star\star)$ we have $0\propto
K\oplus_R P=0\propto I\oplus_R R$. Hence, we conclude that
$0\propto K$ is a finitely generated ideal of $R$. Contradiction,
 since $K$ is
not a finitely generated $D$-module (since $D$ is not a field).
Then $0\propto I$ can not be a strongly Gorenstein projective
ideal of $R$ as desired.\cqfd

\proof\textbf{of Theorem \ref{thmDK}.}  $(1)\Rightarrow (2).$ Let
$I$ be an ideal of $D$ finitely generated strongly Gorenstein
projective. First we show that $I\otimes_D R=I\propto K$ is a
strongly Gorenstein projective ideal of $R$. From
\cite[Proposition 2.12]{BM}, there exists an exact sequence of
$D$-modules: $0\longrightarrow I\longrightarrow P \longrightarrow
I\longrightarrow 0$, where $P$ is finitely generated projective
$D$-module. Moreover, $0\longrightarrow I\propto K\longrightarrow
P\otimes_D R \longrightarrow I\propto K\longrightarrow 0$, is an
exact sequence of $R$-modules. Then,  it remain to prove that
$\Ext_R(I\propto K, R)=0$. From \cite[Proposition 4.1.3]{CE},
$\Ext_R(I\propto K, R)\cong \Ext_D(I, R)$ and  since $R$ is a flat
$D$-module   $\Ext_D(I, R)=0$, \cite[Propositon 2.12]{BM}. Then,
$\Ext_R(I\propto K, R)=0$  and $I\propto K $ is finitely generated
ideal of $R$ strongly Gorenstein projective, so projective. Hence,
$I\cong I\propto K \otimes_R D$ is a projective ideal of
$D$.\\
$(2)\Rightarrow (1).$ Let $J$ be a finitely generated ideal of $R$
strongly Gorenstein projective. It is known from the presentation
of Corollary \ref{corofprop4} that the finitely generated ideals
of $R$ have the forms $I\propto K$  where $I$ is a finitely
generated ideal of $D$ or $0\propto E'$, where $E'$ is a finitely
generated  $D$-submodule of $K$, (without loss of generality we
can assume that $E'$ is a finitely generated ideal of $D$). From
Lemma \ref{lem3.13}, $J=I\propto K$ where $I$ is an ideal finitely
generated of $D$. Then, $J=I\propto K=I\otimes_D R$
 and from \cite[Theorem 2.1]{MO2}, $I$ is a finitely generated ideal
 of $D$ strongly Gorenstein projective, and by hypothesis $I$ is  projective. Then $J=I\otimes_D R$ is a projective ideal of $R$.\cqfd
\end{section}
\begin{section}{Examples}
In this section, we give some examples of rings satisfy the
properties studied in section 2.\\
\indent Recall, from \cite{costa}, a ring $R$ is called an
$(n,d)$-ring if every $R$-module having a finite $n$-presentation
has projective dimension at most $d$. Also from \cite{Mh1}, a
commutative $(n,0)-$ring $R$ is called an $n$-Von Neumann regular
ring. Thus , the $1$-Von Neumann regular rings are the well-know
Von Neumann regular rings. In the following result we show that
the class of $(n,d)$-rings  satisfy the property ``all finitely
generated strongly Gorenstein projective modules are projective".

\begin{thm}\label{propn-RVN}
\begin{enumerate}
    \item Let $R$ be an $(n,d)$-ring. Then, $R$-satisfies
all finitely generated strongly Gorenstein projective $R$-modules
are projective.
    \item Let $R$ be an $n$-Von Neumann regular ring. Then, there is not finitely generated ideal which is strongly   Gorenstein projective.
\end{enumerate}
\end{thm}
\proof
\begin{enumerate}
    \item Let $M$ be a finitely generated strongly Gorenstein projective
$R$-module. From \cite[Proposition 2.12]{BM}, there exists an
exact sequence  of $R$-modules:
$$0\longrightarrow M \longrightarrow
P\longrightarrow M\longrightarrow 0\qquad (\star)$$ where $P$ is
finitely generated projective, then  $M$ is infinitely presented.
Then from \cite[Theorem 2.1]{Mh1}, $\pd_R(M)$ is finite. Thus,
from the exact sequence $(\star)$, we conclude that $M$ is
projective.
    \item  Suppose that $I$ is a strongly Gorenstein projective ideal.
From $(2)$, $I$  is projective, and since $R$ is  $m$-local $I$ is
free. Contradiction since any finitely generated ideal has a
nonzero annihilator (by \cite[Theorem 2.1]{Mh1}).\cqfd
\end{enumerate}

Next, we give examples of rings with  weak global
 dimension  infinite and which satisfies the property ``all finitely
generated strongly Gorenstein projective modules are
 projective".

\begin{cor}\label{corK[[X1]]}
Let $K$ be a field and $A=K[[X_1,X_2,...]]$  with
$m=(X_1,X_2,...)$ the maximal ideal of $A$. And let $R=A\propto
(A/m)^{\infty}$. then:
\begin{enumerate}
    \item $\wdim(R)=\infty$.
    \item $R$ satisfies all finitely generated strongly Gorenstein
    projective $R$-modules are projective.
    \item There is  not finitely generated ideal of $R$ which is
    strongly Gorenstein projective.
\end{enumerate}
\end{cor}
\proof $(1).$ Follows from  \cite[Theorem 3.1]{BMK},  since
$(X_1,X_2,...)$ is a minimal
generating set of $(m)$.\\
$(2)$ and $(3).$ Follows from Theorem \ref{propn-RVN} above and
from \cite[Theorem 2.1]{Mh2}.\cqfd

\begin{cor}\label{thm3.5}
Let $(A,m)$ be an $m$-local ring  such that  $m$ is finitely
generated and $E=(A/m)^{\infty}$. Let $R=A\propto E$ the trivial
ring extension of $A$ by $E$. Then the following conditions holds
in $R$ :
\begin{enumerate}
    \item   $\wdim(R)=\infty$.
    \item  All finitely generated  strongly Gorenstein projective modules are projective.
    \item There is not finitely generated proper ideal which is strongly   Gorenstein projective.
\end{enumerate}
\end{cor}
\proof $(1).$ Follows  from \cite[Theorem 3.1]{BMK}.\\
$(2)$ and $(3).$  Use \cite[Theorem 2.1]{Mh2} and Theorem
\ref{propn-RVN} above.\cqfd

\begin{exmp}\label{exmpRVN}
Let $K$ be a field and $E$  a $K$-vector space of infinite
dimension. Then,  the property,  all finitely generated strongly
Gorenstein projective is projective, holds in $R=K\propto E$.
\end{exmp}

Next, we see an example of ring which satisfies ``all finitely
generated strongly Gorenstein projective ideals are projective.

\begin{exmp}\label{expDK}
Let $(D,m)$ an $m$-local integral domain   and $K=qf(D)$. If
$\wdim(D)$ is finite, then $D\propto K$ satisfies all finitely
generated strongly Gorenstein projective ideals are projective.
\end{exmp}
\proof Follows from Theorem \ref{pro1} and Theorem \ref{thmDK}
above.\cqfd

Next, we see an example of ring $T$ with $\Ggldim(T)=0$, over
which
 the property ``all (finitely generated) strongly
Gorenstein projective modules are projective" does not holds.
Also, this example prove that condition $D$ is not a field in
Theorem \ref{thmDK}  is necessary.

\begin{exmp}\label{contreexp}
Let $K$ be a field and let $T=K\propto K$. Then
    $\Ggldim(T)=0$ but $T$ does not satisfies  all (finitely generated) strongly Gorenstein
projective modules are projective since $0\propto K$ is  a
finitely generated ideal of $T$ strongly
    Gorenstein projective  but it is not
    projective.
\end{exmp}
\end{section}

\end{document}